%\\
%Title: Cycle for integration for zonal spherical function of type $A_n$
%Authors: A.Kazarnovski-Krol
%Comments: 31 pages, Ams-Tex, including 12 Postscript figures
%\\
%Integral of a certain multivalued form  over cycle $\pmb\Delta$ provides
%zonal spherical function of type $A_n$. This paper is devoted to
%quantum group analysis and 
%verification of monodromy properties of the distinguished
%cycle $\pmb\Delta$. Zonal spherical function is a particular conformal
%block of $WA_n$-algebra.
%\\

%This is Ams-Tex

\input amstex.tex
\documentstyle{amsppt}
\input amsppt1.tex
\nologo
\magnification=1000
%\nopagenumbers
%\parindent=0pt
\NoBlackBoxes

%\define\op#1{\operatorname{#1}}

\topmatter
\author   A.Kazarnovski-Krol          \endauthor
\title Cycle for integration for zonal spherical function of type
$A_n$   \endtitle 
\address{
  Department of Mathematics
  Rutgers University
  New Brunswick, NJ 08854, USA}
\abstract{
Integral of a certain multivalued form  over cycle $\pmb\Delta$ provides
zonal spherical function of type $A_n$. This paper is devoted to
quantum group analysis and 
verification of monodromy properties of the distinguished
cycle $\pmb\Delta$. Zonal spherical function is a particular conformal
block of $WA_n$-algebra
 }
\endtopmatter

\document

\subhead 0.0 Notations \endsubhead

$\alpha_1, \alpha_2, \ldots, \alpha_n$ - simple roots of root system
of type $A_n$

$\Cal R_{+}$ - set of positive roots

$\Cal R$ root system of type $A_n$

$\delta= {1\over 2} \sum_{\alpha \in \Cal R_{+}} \alpha$ -halfsum of
positive roots

$k$- complex parameter ( `halfmultiplicity' of a root)
$$\rho = {k \over 2} \sum_{\alpha \in \Cal R_{+}} \alpha$$

$\varkappa= - {1\over k}$

$q= exp( {{2 \pi i} \over{\varkappa}})$

$\eta_1, \eta_2, \ldots, \eta_n$ are fundamental weights: $(\eta_i,
\alpha_j)=\delta_{ij}$

$h_1, h_2, \ldots, h_{n+1}$ are weights of the vector representation
of $sl(n+1)$

$h_1=\eta_1$

$h_2=h_1-\alpha_1, h_3=h_1-\alpha_1-\alpha_2, \ldots,
h_{n+1}=h_1-\alpha_1-\ldots -\alpha_n$

$$h_1 +h_2 +\ldots + h_{n+1}=0  $$

$R$- R-matrix

$\pmb\Delta$ - distinguished cycle for integration for zonal spherical
function of type $A_n$

$\Delta$ - comultiplication in quantum group

\head 0. Introduction\endhead

Let $G$ be a  connected real semisimple Lie group with  finite center, 
K -its maximal compact subgroup,
$G/K$- Riemannian symmetric space. Let $T_g=T_g^{\lambda}, g \in G$  be a
continuous unitary representation of $G$ acting in a Hilbert space $H$
, which contains a spherical
vector $\xi$, i.e. $K \xi =\xi$ and assume that
$(\xi,\xi)=1$. $\lambda$ is a parameter defining this representation.
Then the function $\phi_{\lambda}(g)=({T_g}^{\lambda} \xi,\xi)$ is called zonal
spherical function [3]. In particular, $\phi_{\lambda}(e)=1$ and it is
right and left $K$-invariant.
Zonal spherical function is a common eigenfunction of Laplace-Casimir
operators :
$$\Cal L \phi_{\lambda}(g)= \gamma(\Cal L)( i \lambda) \phi_{\lambda}(g)$$
where $\gamma(\Cal L)( i \lambda)$ is a homomorphism of Laplace -Casimir
operators into complex numbers.
Using Cartan decomposition $G=KAK$ zonal spherical function
$\phi_{\lambda}(g)$ is considered as a function on $A$.
Let also $\frak a=Lie(A)$, then $\lambda \in \frak a^{*}$.

Restriction of zonal spherical function to $A$ is a common
eigenfunction of radial parts $\overset \circ \to {\Cal  L}$ of Laplace-Casimir
operators $\Cal L$:
$$\overset \circ \to {\Cal L} (\phi_{\lambda}(a))= \gamma(\overset \circ
\to {\Cal L})(i \lambda)\phi_{\lambda}(a)           \tag{1}$$
where $\gamma $ is Harish-Chandra homomorphism, $a \in A$.
Among others the operator of second order plays the predominant 
role:
$$\overset \circ \to {\Cal L_2}=
H_1^2+\ldots +H_l^2 
+ \sum_{\alpha \in \Cal R_{+}} m_{\alpha} \frac{e^{\alpha}+e^{-\alpha}}
{e^{\alpha}-e^{-\alpha}}H_{\alpha}$$
where $m_{\alpha}$ is a multiplicity of restricted root $\alpha$ [6].
It turns out that one can consider $m_{\alpha}$ as parameters, but the
condition $m_{\alpha}=m_{w \alpha}$ is required. So there are as many
independent parameters as orbits of the Weyl group in restricted root system.
One may start with second order differential operator with generic
parameters,
then recover
the whole system of differential operators (radial parts)[41,20].
This system turns out to be holonomic, locally it has
$|W|$-dimensional space of solutions,
where $|W|$ is the cardinality of the Weyl group $W$,
cf. corollary 3.9 of [20 ].
Among those solutions there is a distinguished one which corresponds
to zonal spherical function. It is characterized by analyticity at
unity
cf. theorem 6.9 [20]. 
  
We restrict ourselves to the case of root sytem of type $A_n$. 
For the second order differential operator we use
$$
\overset \circ \to {\Cal L_2}=
\sum_{i} \bigg(  z_i \frac{\partial} {\partial
z_i}\bigg)^2 - k \sum_{i<j} \frac{z_j +z_i}{z_j-z_i}\bigg(z_i
\frac{\partial}{\partial z_i}- z_j\frac{\partial}{\partial z_j}\bigg) \quad. 
$$
In [17 ] we provided an integral representation for the solutions of
system (1), in [18 ] we described  cycle $\pmb\Delta$ for integration
for zonal spherical function and obtained an explicit version of
Harish-Chandra decomposition. 
\smallskip
This  paper is devoted to quantum group analysis and verification of
monodromy properties of distinguished cycle $\pmb\Delta$ 
(cycle $\pmb\Delta$ is recalled in 2.1 below).  
Cycle  $\pmb\Delta$ is a contour for integration for zonal spherical
function  of type $A_n$ of a suitable multivalued form.
The form is of type considered in refs. [8,11] and thus all the
machinery including quantum groups and R-matrices can be applied.
We obtained the form using very simple principle:
there is only one line which passes through two given points.
Zonal spherical functions for $SL(n,\Bbb C)$ are calculated by I. Gelfand
and M. Naimark   and using the same method can be calculated for 
$SL(n,\Bbb R)$ (see  [1]). This provides us with the form for
parameters
 $k=1$ and 
$k=\frac{1}{2}$, where $k$ is a halfmultiplicity of restricted root.
Now use the principle and extend powers of factors  linearly on $k$ .
The obtained form has several advantages:
  cycle for integration for zonal
spherical function is real and compact; the number of variables of
integration is independent on parameter $\lambda$ (appearance of the
flag  manifold); there is no complicated meromorphic
factor.
\smallskip  
Note: in refs. [15,54] the  system $(1)$ is proved to be related 
to a particular  case of trigonometric Knizhnik-Zamolodchikov
equation. In particular, this implies that solutions to the system $(1)$ 
can be obtained from the solutions of Knizhnik-Zamolodchikov equations
by symmetrization procedure. Solutions to Knizhnik-Zamolodchikov
equations are given by certain multidimensional integrals, whose
integrand has the standard part times complicated meromorphic factor.
This complicated meromorphic factor becomes even more complicated
after symmetrization. We would like to emphasize that in this
particular case this unpleasent meromorphic factor is not needed,
 cf. [17] , theorem 6.3. 

\smallskip
Knizhnik-Zamolodchikov equations are originated in WZNW theory.
Reduction from WZNW to $WA_n$-algebras is well
discussed in the literature [75,63] ( quantum Drinfeld-Sokolov
reduction).
Results of [17] (absence of meromorphic factor)
imply that solutions to the system $(1)$ (which is isomorphic to 
Calogero-Sutherland model)  
are provided by the conformal blocks of $WA_n$ -algebra.
\smallskip
Quantum group approach assumes the following. With the multivalued
form one associates tensor product of Verma modules over quantum
group . Homology of certain type of discriminantal configuration
[8,71] are described by the singular vectors
of the tensor product of irreducible highest weight modules. 
 Half-monodromy(=braiding) is given by the 
R-matrix (PR, where P is a permutation). Universal R-matrix is
provided by  the Drinfeld's double [26]. 
\smallskip
Here is the organization of the paper.
In sections 1-3 we recall the multivalued form, distinguished cycle
$\pmb\Delta$, and normalization constant of ref. [16].
In section 4 we recall the version of quantum group used  in refs. [8, 11]
for the explicit version of Kohno's theorem:
half-monodromy=R-matrix [28]. In section 5 we encode the distinguished
cycle $\pmb\Delta$ as an element of the corresponding tensor product 
and check the monodromy properties.
 Cycle $\pmb\Delta$ has the meaning of q-antisymmetric tensor, cf.
theorem 5.7 below, and corresponds to a particular conformal block,
fig. 10.
Finally, we discuss the properties of the tensor
product with vector representation.
\smallskip
\head 1. Multivalued form and discriminantal configuration \endhead

Consider the following set of variables:

$z_l, \;\;l=1, \ldots, n+1, \; t_{ij} ,\; i=1, \ldots, j,\;\;
j=1,\ldots, n.$

Variables $z_l$ have meaning of arguments, while variables
$t_{ij}$ are variables of integration.

It is convenient to organize variables $z_l,\; t_{ij}$ in the
 form of a pattern, cf. fig 1. The idea of such an organization is
borrowed from Gelfand-Zetlin patterns [7].

\midinsert
$$
\matrix
z_{1}&&z_{2}&&\ldots&&\ldots&&z_{n+1}\\\\
&t_{1,n}&&t_{2,n}&&\ldots&&t_{n,n}&\\\\
&&\ldots&&\ldots&&\ldots&&\\\\
&&&t_{1,2}&&t_{2,2}&&&\\\\
&&&&t_{1,1}&&&&
\endmatrix
$$
\botcaption{Figure 1}
Variables organized in a pattern 
\endcaption
\endinsert

\definition{Definition 1.1} Consider the following multivalued form
$\omega(z,t):$

$$
\align
\omega(z,t):= &\prod_{i=1}^{n+1} z_i^{\lambda_1 +{k n\over 2}} 
 \prod_{i_1 > i_2}{(z_{i_1}-z_{i_2})^{1-2k}}\\
 &\times\prod_{l=1}^{n+1}\prod_{i=1}^n{(z_l-t_{i,n})}^{k-1}\\
 &\times\prod_{j=1}^{n-1} 
 \prod_{i_1=1}^{j+1} \prod_{i=1}^j ( t_{ij} - t_{i_1,j+1} )^{k-1}\\
 &\times\prod_{j=2}^n \prod_{i_1>i_2} {(t_{i_1,j}-t_{i_2,j})^{2-2k} }\\
 &\times\prod_{j=1}^n \prod_{i=1}^{j} {t_{ij}^{\lambda_{n-j+2}-
 \lambda_{n-j+1} -k} } \quad
 {dt_{11} dt_{12} dt_{22} \ldots dt_{nn} }
\endalign
$$

Here  $k$ is a complex parameter - `halfmultiplicity' of a 
root, $\lambda_1, \ldots, \lambda_{n+1}$ are complex parameters subject to
the homogeneity condition:
$$\lambda_1 +\lambda_2 +\ldots +\lambda_{n+1}=0$$

\enddefinition

Before proceeding further we would like to make a convention.

\subhead{Convention 1.2}\endsubhead
A complex number $z$ can be represented as
$z=r e^{i \alpha}$, where $r, \alpha$ are real numbers, $r \ge 0$. 
$r$ is called
absolute value of $z$, while $\alpha$ is called the phase of $z$. When we say 
that the phase of a complex number $z$ is equal to $0$, we mean
that $\alpha=0$, or the number itself is real and nonnegative. 
\smallskip

\subhead{1.3 Configuration }\endsubhead

Let  $m={n(n+1) \over 2}$. Consider $(n+1+m)$-dimensional complex  space 
$\Bbb C^{n+1+m}$
with coordinates 
$z_1,z_2,\ldots,z_{n+1}, t_{11} , t_{12}, t_{22}, \ldots , t_{nn}$.
Let's delete the following hyperplanes:
$$ t_{i_1,j}-t_{i_2,j}=0 \quad i_1 <i_2,\; j=1, \ldots, n$$
$$t_{i_1,j}- t_{i_2, j+1}=0  \quad j=1, \ldots, n-1$$
$$z_{i_1}-t_{i_2,n}=0 \quad i_1=1,\ldots, n+1 ;\; i_2=1,\ldots,n $$
$$t_{ij}=0 \quad i=1,\ldots, j; \;j=1,\ldots, n $$
$$z_i-z_j=0 \quad i < j $$
$$z_i=0 \quad i=1,\ldots,n+1.$$

Denote the complement of $\Bbb C^{m+n+1}$ to the union of  above 
hyperplanes by
$U_{n+1+m}$.

Denote by $Loc_{\mu}  $ the trivial 1-dimensional bundle over
$U_{n+1+m}$
with the integrable 
connection  $\nabla_{\mu}$ with the connection form
$$
\align
& \sum (k-1) \frac{d(t_{jn}-z_i)}{t_{jn}-z_i} +
\sum (1-2k) \frac{d (z_i - z_j)}{z_i-z_j}+\sum (\lambda_1 +
\frac{kn}{2})\frac{d z_i}{z_i}+
\sum (2-2k) \frac{d (t_{i_1,j}-t_{i_2,j})}{t_{i_1,j}-t_{i_2,j}}   + \\
& \sum (k-1)\frac{d(t_{i_1,j} -t_{i_2,j+1})}{t_{i_1,j}-t_{i_2,j+1}}  +
\sum  (\lambda_{n+2-j}-\lambda_{n+1-j}-k) \frac {d t_{ij}}{t_{ij}} 
\endalign
$$
Denote by $S_{\mu}$ the local system of horizontal sections of
$\nabla_{\mu}$.
Consider the projection on the first $n+1$ coordinates 
$\Bbb C^{n+1+m} \mapsto \Bbb C^{n+1}$.
For $z=(z_1,z_2,\ldots,z_{n+1})$ such that $z_i \ne z_j$ for all 
$i < j$ set 
$$U(z)=\Big\{ {(\tilde{z},t) \in U_{n+1+m}} | \tilde{z} = z \Big\} \; .$$
Restrictions of $Loc_{\mu}, S_{\mu}$ to $U(z)$ are denoted
by $Loc_{\mu}(z), S_{\mu}(z)$. 
Denote by $S^{*}$ the dual local system and consider homology of $U(z)$
 with coefficients in $S^*$  (extended by $!$ as explained in ref. [8]).
Configuration is preserved under the action of the product of symmetric
groups :
$$\Sigma =S_n \times S_{n-1} \times\ldots \times S_2 , $$  
where $S_j$ permutes $ t_{1j}, t_{2j}, \ldots, t_{jj}$.
 Then one considers the antiinvariant part of the homology
group  with respect to the action of $\Sigma$: \;  $H_{*!,m}(U,S^{*})^{-}$.
\smallskip
Remarkably, in order to calculate the cohomology group of the local
system of the complement to finite set of hyperlanes in the
nonresonance case one can use
finite-dimensional complex of hypergeometric forms in the spirit of 
Arnold-Orlik-Solomon, cf. [12  ,67].

\head 2. The distinguished cycle $\pmb\Delta$ \endhead

Assume that $z_1, z_2, \ldots , z_{n+1}$ are real and 

$$0 < z_1 < z_2 < \ldots < z_{n+1}.$$

\definition{Definition 2.1} Define cycle $\pmb\Delta=
\pmb\Delta(z)$ by the following inequalities:

$t_{i,j+1} \le t_{ij} \le t_{i+1, j+1}$ and 

$z_i \le t_{in} \le z_{i+1}$  cf. [16], definition 2.1.

Define form $\omega_{\Delta}(z,t)$ as:

$$
\align
\omega_{\Delta}(z,t):= &\prod_{i=1}^{n+1} z_i^{\lambda_1 +{k n\over 2}} 
 \prod_{i_1 > i_2}{(z_{i_1}-z_{i_2})^{1-2k}}\\
 &\times\prod_{i \le l}{(z_l-t_{i,n})}^{k-1}\prod_{i > l}{(t_{i,n}-z_l)}^{k-1}  \\
 &\times\prod_{j=1}^{n-1} 
 \prod_{i_1 > i_2} ( t_{i_1,j} - t_{i_2,j+1} )^{k-1} 
 \prod_{i_2 \ge  i_1} ( t_{i_{2},j+1 } - t_{i_1,j} )^{k-1}\\
 &\times\prod_{j=2}^n \prod_{i_1>i_2} {(t_{i_1,j}-t_{i_2,j})^{2-2k} }\\
 &\times\prod_{j=1}^n \prod_{i=1}^{j} {t_{ij}^{\lambda_{n-j+2}-
 \lambda_{n-j+1} -k} } \quad
 {dt_{11} dt_{12} dt_{22} \ldots dt_{nn} }
\endalign
$$

It is assumed that phases of factors in the formula for
$\omega_{\Delta}$
are equal to zero if $k$ and $\lambda_1, \lambda_2, \ldots,
\lambda_{n+1} $ are real. In other words we choose the section of the
local system to be positive over $\pmb\Delta$ if $\lambda, k$ are
real. this geometric definition of cycle $\pmb\Delta$ is justified 
by theorem 5.7 below, i.e. $\pmb\Delta$ is really an element of
$H_{*!,m}(U,S^{*})^{-}$. This geometric definition is motivated by the 
classical calculations of Gelfand and Naimark of zonal spherical
function for $SL(n,\Bbb C)$. 
 
\enddefinition

\smallskip
\head{3. Analytic considerations} \endhead

Let  
$${\lambda_{n-j+2}} -{\lambda_{n-j+1}} -k =0$$
  for $j=1, \ldots , n$ and 
$\lambda_1+{k n\over 2}=0$, i.e we kill an affine part.
% or ,in other
%words, we erase the first factor in the tensor product of Verma
%modules.

Then in these hypotheses 
$$
\int_{{\pmb\Delta(z)}} \omega_{\Delta}(z,t) =
{
 {\Gamma(k)\Gamma(k)^2\ldots\Gamma(k)^{n+1}}
 \over
 {\Gamma(k)\Gamma(2k) \ldots\Gamma((n+1) k)}
}
$$
cf. ref. [16] theorem 1.5  and remark 1.6.

  Following the classical work of I.M. Gelfand and 
M.A. Naimark cf. [23]
 ,  let
$$
\tau_{ij}=
{
 {\prod\limits_{i_1=1}^{j-1}(t_{i_1,j-1}-t_{ij})}
 \over
 {\prod\limits_{i_1\ne i}(t_{i_1,j}-t_{ij})}
}\;\; ,
$$
$i=1,\ldots,j \;$ ;
$j=2,\ldots,n$.

Note that ${\sum\limits_{i=1}^j}\tau_{ij}=1$, and
$$
{
 {D(\tau_{1,j},\ldots,\tau_{j-1,j})}
 \over
 {D(t_{1,j-1},\ldots,t_{j-1,j-1})}
}
=
{
 {{\prod\limits_{1\le i<k\le j-1}}(t_{i,j-1}-t_{k,j-1})}
 \over
 {{\prod\limits_{1\le i<p\le j}}(t_{ij}-t_{pj})}
}
$$
[see \cite{23} for the details].
Let also 
$$
\tau_{i,n+1}={
           {{\prod\limits_{i=1}^{n}}(t_{i_1,n}-z_i)}
           \over
           {{\prod\limits_{i_1\ne i}}(z_{i_1}-z_i)}
          }
\;\;\;\;\;\;\;\;\;\;\;  i=1,\ldots,n+1.
$$

One has ${\sum\limits_{i=1}^{n+1}}\tau_{i,n+1}=1$ and
$$
{
 {D(\tau_{1,n+1},\ldots,\tau_{n,n+1})}
 \over
 {D(t_{1,n},\ldots,t_{n,n})}
}
=
{
 {{\prod\limits_{1\le i<k \le n}}(t_{i,n-1}-t_{k,n})}
 \over
 {{\prod\limits_{1\le i<p\le n+1}}(z_i-z_p)}
}
$$
In  variables $\tau_{ij}\;\;\; i=1,\ldots,j-1,\;\; j=2,\ldots,n+1$
integral  $ \int_{\pmb\Delta} \omega_{\Delta}$ is written as:
$$
\align
\int_{\Delta} \omega_{\Delta}
=
&\int
{
{\prod\limits_{j=1}^{n+1}} 
 (
  (\tau_{1j}\tau_{2j}\ldots\tau_{j-1,j})
  (1-\tau_{1j}-\ldots -\tau_{j-1,j})
 )^{k-1}
}\\
& \times d\tau_{12}d\tau_{13}d\tau_{23}\ldots d\tau_{1,n+1}\ldots
d\tau_{n,n+1}\;\; .
\endalign
$$

Remarkably in variables $\tau_{ij}$ the integration is being performed
over one-dimensional simplex times two-dimensional simplex times so
on, times $n$-dimensional simplex. Here one-dimensional simplex
corresponds to a line in two-dimensional plane, two-dimensional simplex 
corresponds to two-dimensional plane in three-dimensional plane and so
on. Thus the integral remembers the flag manifold.

So using Dirichlet's formula one gets
$$
\int_{\Delta}\omega_{\Delta}=
{
 {\Gamma(k)\Gamma(k)^2\ldots\Gamma(k)^{n+1}}
 \over
 {\Gamma(k)\Gamma(2k) \ldots\Gamma((n+1)k)}
}
$$

 The constant does not depend on $z_i$ at all and surely
remains the same under analytic continuation.
This is  nontrivial since  form $\omega=\omega(z,t)$
and cycle $\pmb\Delta={\pmb\Delta}(z)$ do depend on 
$z=(z_1,z_2,\ldots,z_{n+1})$.

\remark{Remark 3.1} In view of sections 4 and 5  below, killing of affine part 
corresponds to erasing of the first factor in the tensor product 
of Verma modules.
\endremark
\medskip
\head 4. Quantum group \endhead

Quantum groups are introduced by Drinfeld [26], Jimbo [27], Kulish,
Reshetikhin, Sklyanin [37]. 
We briefly  recall the necessary material from [8,11] and refer
directly to these references for more details. See also [29]. 

\subhead {4.1 Root system  }\endsubhead

Let $\Bbb R^{n+2}$ be Euclidean $(n+2)$-dimensional vector space
 with inner
product $(.,.)$ and 
with $g_0,g_1,\ldots,g_{n+1}$ as the orthonormal basis  .
Let's realize simple roots of root system of type $A_n$ 
as
$\alpha_i=g_i-g_{i+1}$ for $i=1,\ldots,n$.
Set also $\alpha_0= g_0-g_1 $.

In particular, one has
$$(\alpha_i, \alpha_i)=2$$
$$(\alpha_i, \alpha_j)=0 \quad \text{for} \quad |i-j| >1$$
$$(\alpha_i, \alpha_j)=-1 \quad \text{for} \quad |i-j|=1$$
\smallskip 
Set also
$$\alpha^{\vee}=\frac{2 \alpha}{(\alpha,\alpha)}.$$

\subhead{ 4.2. Quantum group $U_q(sl(n+2))$} \endsubhead

Consider $\Bbb C $-algebra with generators $e_i, f_i,$ and
$K_i^{1\over 2}, K_i^{-1\over 2}$ 
,subject to the relations:
$$
 K_j^{1\over 2} e_i=q^{\frac{(\alpha_i,\alpha_j^{\vee})}{4}}e_i K_j^{1\over 2} \quad 
$$

$$
K_j^{1\over 2} f_i=q^{-\frac{(\alpha_i,\alpha_j^{\vee})}{4}} f_i K_j^{1\over 2} 
$$

$$[e_i,f_j]=(K_i-K_i^{-1}) \delta_{ij}$$

$$K_i^{{\pm} {1\over 2}} K_j^{1\over 2} =K_j^{1\over 2} K_i^{{\pm}
{1\over 2}}$$

$$K_i^{1\over 2} K_i^{-1\over 2}=K_i^{-1\over 2} K_i^{1\over 2} = 1$$ 
Comultiplication is defined by the rule
$$
\Delta(K_i^{{\pm} {1\over 2}})=K_i^{{\pm} {1\over 2}} 
\otimes K_i^{{\pm} {1\over 2}}$$
$$ \Delta(f_i)=f_i\otimes K_i^{1\over 2} + K_i^{-{1\over 2}} \otimes f_i
$$
$$ \Delta(e_i)=e_i \otimes K_i^{1\over 2} + 
K_i^{-{1\over 2}} \otimes e_i$$
\smallskip

The following are the quantum Serre's relations:

$$
f_i^2 f_{i+1} -(q^{\frac{1} {2}}+ q^{-\frac{1}{2}})f_i f_{i+1} f_i
+f_{i+1}f_i^2=0
$$

$$
f_{i+1}^2 f_i -( q^{\frac{1}{2}}+ q^{-\frac{1}{2}})f_{i+1} f_i
f_{i+1}+ f_i f_{i+1}^2=0
$$

$$f_i f_j =f_i f_j \quad \text{for} \quad |i-j| \ne 1$$

$$
e_i^2 e_{i+1} -(q^{\frac{1}{2}}+q^{-\frac{1}{2}})e_i e_{i+1} e_i
+e_{i+1} e_i^2=0
$$

$$
e_{i+1}^2 e_i -(q^{\frac{1}{2}} +q^{-\frac{1}{2}})e_{i+1} e_i e_{i+1}
+e_i e_{i+1}^2=0
$$

$$e_i e_j =e_j e_i \quad \text{for} \quad |i-j| \ne 1$$

$\Bbb C$-algebra generated by elements $e_i,f_i,K_i^{1\over 2},
K_i^{-{1\over 2}} $ 
subject to the above relations and with the comultiplication will be
referred to as quantum group $U_q(sl(n+2))$. Antipode $S$ and counit
$\epsilon$ are  defined appropriately.

\smallskip
\subhead{4.3. Verma module} \endsubhead
For $\Lambda \in span\{\alpha_i\}$ denote by $M(\Lambda)$
the Verma module generated over quantum group by a single vector 
$v$ subject to relations: 
$K_i^{ 1\over 2} v= q^{\frac{(\Lambda,\alpha_i^{\vee})}{4}}v $, 
$e_i v=0$ for all $i$.

$$
M(\Lambda)_{\mu}=\Big\{x \in M(\Lambda) \;| K_i^{1\over 2} x=
q^{\frac{(\Lambda-\mu,\alpha_i^{\vee})}{4}}x \Big\}.
$$

Set $\tau(e_i)=f_i, \; \tau(f_i)=e_i, \; \tau(K_i^{{\pm}{1\over 2}})
=K_i^{{\pm}{1\over 2}}$ on generators and
extend $\tau$ as algebra antihomomorphism.

Put $M(\Lambda)^*= \oplus_{\mu} {M(\Lambda)_{\mu}}^*$.
Define the structure of  quantum group module on $M(\Lambda)^*$ by the rule
$(g \phi,x)=(\phi, \tau(g) x)$.

\subhead{4.4. Contravariant form}\endsubhead
Contravariant form on a Verma module with highest weight
vector $v$ is defined such that

$$S(v,v)=1$$
$$S(f_i x,y)=S(x,e_i y)$$
for all $i, x, y$.
\smallskip
The form $S$ defines the homomorphism of modules:
$$S \: M(\Lambda) \mapsto M(\Lambda)^*.$$

\smallskip
Example 1.
$$S \: v \mapsto v^*.$$

Example 2. 
$$
S \: f_1 v \mapsto
(q^{\frac{(\Lambda, \alpha_1)}{2}}-q^{-\frac{(\Lambda, \alpha_1)}{2}}) (f_1v)^*
,$$
see also fig. 2.
\smallskip
Let $$L(\Lambda)= M(\Lambda)/ {\text{Ker} S}$$ be irreducible 
module with highest weight $\Lambda$.
\medskip
\subhead{4.5. R-matrix}\endsubhead
R-matrix is defined by the following expression:
$$
R=\sum_{\mu} q^{\frac{\Omega_0}{2} +\frac{1}{4}({\mu} \otimes 1-
1 \otimes {\mu}) +d(\mu)}\Omega_{\mu} 
$$
cf. [26].
Here $ \Omega_0 $ is the element
corresponding to the inner product $(.,.)$ ; for
$$
\mu =
l_0{{\alpha}_0}+{l_1}  {{\alpha}_1} +\ldots + {l_n} {{\alpha}_n}  
,$$
where $l_0, l_1, \ldots, l_n$ are nonnegative integers, 
$ d(\mu) \in \Bbb C $ is a constant defined as follows:
represent $\mu=\sum {l_i} {{\alpha}_i} $ as a sum of simple roots 
with repetitions
$ \mu ={{\alpha}_{i_1}} +{{\alpha}_{i_2}} +\ldots + {{\alpha}_{i_n}} $.
Then 
$$ d(\mu)= -\sum_{p \le q} \frac{(\alpha_{i_p},\alpha_{i_q})}{4} .$$
$ \Omega_{\mu} $ is a canonical element cf. [26].

 $R$ defines a linear operator 
$$R \: M\otimes M^{\prime} \mapsto M \otimes M^{\prime}.$$
The following diagram is commutative :
$$
\CD
M\otimes M^{\prime} @ >R >> M \otimes M^{\prime}\\
@V S VV@VVS V\\
M^*\otimes M^{\prime *}@>> R^* > M^* \otimes M^{\prime *}\\
\endCD
$$
 cf. theorem 7.6.8 of ref.[11].
$R$  induces a homomorphism of irreducible highest weight modules:
$$R \: L(\Lambda(1)) \otimes L(\Lambda(2)) \mapsto L(\Lambda(1)) \otimes 
L(\Lambda(2))$$
which will be also denoted by $R$.

Denote by $P$ the transposition of two factors in the tensor product:
$$P \: M \otimes M^{\prime} \mapsto M^{\prime} \otimes M .$$

\head 5. Quantum group and cycle $\pmb\Delta$ \endhead

\subhead 5.1 Data \endsubhead
  We are going to check  the monodromy
properties of cycle  $\pmb\Delta$ using quantum group argument,
cf. [8, 11,29].
\smallskip
Take now a different indexation of variables $t_{ij}$. Namely, we are
going to use
 $ \{ {t_i}^{(j)} \;| j=1,\ldots,n ;\; i=j,\ldots ,n \} $ cf.[15].
Set also $z_0=0$ ( affine configuration).
Consider the following multivalued form:
\smallskip
$$
\multline
\Omega(z,t)= \prod(z_0-z_i)^{(\Lambda(0), \Lambda(i))\over \varkappa}
\prod (z_i -z_j)^{(\Lambda(i),\Lambda(j))\over \varkappa} \\
\times \prod (z_i -{t_l}^{(j)})^{(\Lambda(i), -\alpha_j)\over \varkappa}
\prod (z_0- {t_l}^{(j)})^{(\Lambda(0), -\alpha_j) \over \varkappa}
\prod (t_l^{(j)}-t_{l^\prime}^{(j^\prime)})^{(-\alpha_j,
-\alpha_{j^\prime})\over \varkappa} d{t_1}^{(1)} \ldots d{t_n}^{(n)}
\endmultline
$$
\smallskip
Integrals of forms of this type are considered in refs. [8,11,...].
Now we would like to specialize $\Lambda(0),\Lambda(1), \ldots,
\Lambda({n+1}), \alpha_1, \alpha_2, \ldots, \alpha_n$ as follows.

Recall that $\Bbb R^{n+2}$ is an $(n+2)$-dimensional Euclidean vector
space 
with
$g_0,g_1, \ldots, g_{n+1}$ as the orthonormal basis.
For $i=1,\dots,n+1$ set 
$$\Lambda(i) = \Lambda = g_1 - g_0$$ i.e. to each variable $z_i, i=1,\ldots,n+1$
assign the same vector $g_1-g_0$.

Recall that simple roots of root system of type $ A_n$ are realized as
follows:
$$\alpha_i= g_i - g_{i+1} \;,\quad \text{for} \quad i=1,\ldots,n.$$
\medskip
\remark{ Remark 5.2 } Note: the projection of $\Lambda(i)$ on the 
span of $\alpha_i$, $i=1,\ldots ,n$ is exactly the first fundamental
weight, i.e. $(\Lambda(i), \alpha_j)=\delta_{1j}$, where $\delta_{1j}$
is a Kronecker's delta, but $(\Lambda(i),\Lambda(i))=2$. This is not
very important since it changes only the power of $\prod (z_i-z_j)$
before the integral.  
\endremark
\smallskip
Let $\Cal R$ denotes the root sytem of type $A_n$ with simple roots
$\alpha_1, \alpha_2, \ldots, \alpha_n$ as before.
$\Cal R_{+}$ denotes the set of positive roots.
Let $\delta$ be half the sum of positive roots:
$$\delta = \frac{1}{2} \sum_{\alpha \in \Cal R_{+}} \alpha$$
$$\rho = {k \over 2} \sum_{\alpha \in \Cal R_{+}} \alpha$$ 
Set also $\varkappa = -{1\over k}$.
Let $\lambda$ belongs to the span of ${\alpha_1,\ldots , \alpha_n}$.
Set $\Lambda(0)=\varkappa \lambda- \delta$, so that
$${\Lambda(0) \over \varkappa} = \lambda + \rho$$
Let also $$\overline{\lambda}= \varkappa \lambda$$
\smallskip
 The multivalued form $\omega(z,t)$ of  section 1 differs from $\Omega(z,t)$
in the above setting only by some meromorphic factor  
which does not contribute to the
monodromy and thus can be omitted for the purposes of this section.
\smallskip
Set $q= exp(\frac {2 \pi i }{\varkappa})$.
\smallskip
\remark{Remark 5.3}Note: our form is slightly different from the
form
 of ref. [15],
in particular, in setting of  [15] cycle $\pmb\Delta$ will not serve
as a cycle for integration for zonal spherical function (also
 we do not have the complicated meromorphic factor).
For example, 
$$
\int_{z_1}^{z_2} (t-z_1)^{-k}(z_2-t)^{-k}dt =
(z_2-z_1)^{1-2 k} \quad  \frac{\Gamma(1-k) \Gamma(1-k)}{\Gamma(2-2k)}
$$
In particular, if $z_2$ goes around $z_1$ counterclockwise
then  it earns the factor  $exp(2 \pi i (-2 k))$ and thus the cycle is not
preserved under the monodromy.
\endremark
\smallskip
\subhead {5.4. The most trivial example }\endsubhead
Before proceeding furhter we want to consider the most trivial
example. Namely , consider 
$$(z_1 -z_2)^{(\Lambda(1),\Lambda(2)) \over \varkappa}$$
If $z_2$ goes around $z_1$ counterclockwise then this function 
earns the factor  $exp( {2 \pi i \over \varkappa}(\Lambda(1),\Lambda(2)))$.
If $z_2$ goes halfway around $z_1$ then the function earns the factor
$exp( {\pi i \over \varkappa} (\Lambda(1), \Lambda(2)))=q^{(\Lambda(1),
\Lambda(2))\over 2}$.
At the same time  consider $v_1 \otimes v_2$ tensor product of highest
weight vectors
of
modules of the corresponding quantum group
of weights $\Lambda(1)$ and $\Lambda(2)$ correspondingly.
Let $R$ be the $R$-matrix, then 
$$R(v_1 \otimes v_2)= q^{\Omega_0 \over 2} v_1 \otimes v_2,$$
where $\Omega_0$ is the canonical element corresponding to inner
product
,i.e. 
$$q^{\Omega_0 \over 2} v_1 \otimes v_2 = q^{(\Lambda(1),
\Lambda(2))\over 2} v_1 \otimes v_2 $$
in agreement with the above considerations.
\smallskip
\subhead 5.5. Case of root system of type $A_1$ $(n=1)$ \endsubhead

Consider the tensor product of three dual Verma modules over $ U_q(sl(3))$
(simple roots $\alpha_0, \alpha_1$) 
$$ M(\Lambda(0))^{*} \otimes M(\Lambda(1))^{*} \otimes M(\Lambda(2))^{*}.$$

Let $v_0 \otimes v_1 \otimes v_2$ be the tensor product of highest
weight vectors.
Then cycle $\pmb\Delta$ is encoded as:
$$
\multline
v_{\Delta}= - q^{(\Lambda(1), -\alpha_1) \over 4} v_0^{*} \otimes (f_1 v_1)^{*}
\otimes v_2^{*} 
+ q^{-{(\Lambda(2), -\alpha_1) \over 4}} v_0^{*} \otimes v_1^{*} \otimes
(f_1 v_2)^{*}= \\
- q^{-{1\over 4}} v_0^{*} \otimes (f_1 v_1)^{*} \otimes v_2^{*} +
q^{ 1 \over 4} v_0^{*} \otimes  v_1^{*} \otimes (f_1 v_2)^{*} .
\endmultline
$$
Here * means dual with respect to the contravariant  form $S$.

Consider the action of R-matrix on the second and third component of
tensor product 
$$
R \: M(\Lambda(1))^{*} \otimes M(\Lambda(2))^{*} \mapsto  M(\Lambda(1))^{*} \otimes
M(\Lambda(2)) ^{*}.
$$
Recall that $P$ denotes permutation of factors:
$$
P \: M(\Lambda(1))^{*} \otimes M(\Lambda(2))^{*} \mapsto M(\Lambda(2))^{*} 
\otimes M(\Lambda(1))^{*}  .$$
 
Now one can utilize formulas of example 1 of [8].

Then 
$$P R \: v_{\Delta} \mapsto (-1)v_{\Delta}$$
And so
$$
(P R)^2 \: v_{\Delta}  \mapsto v_{\Delta}.$$

\midinsert\vskip 0.5cm
\includegraphics{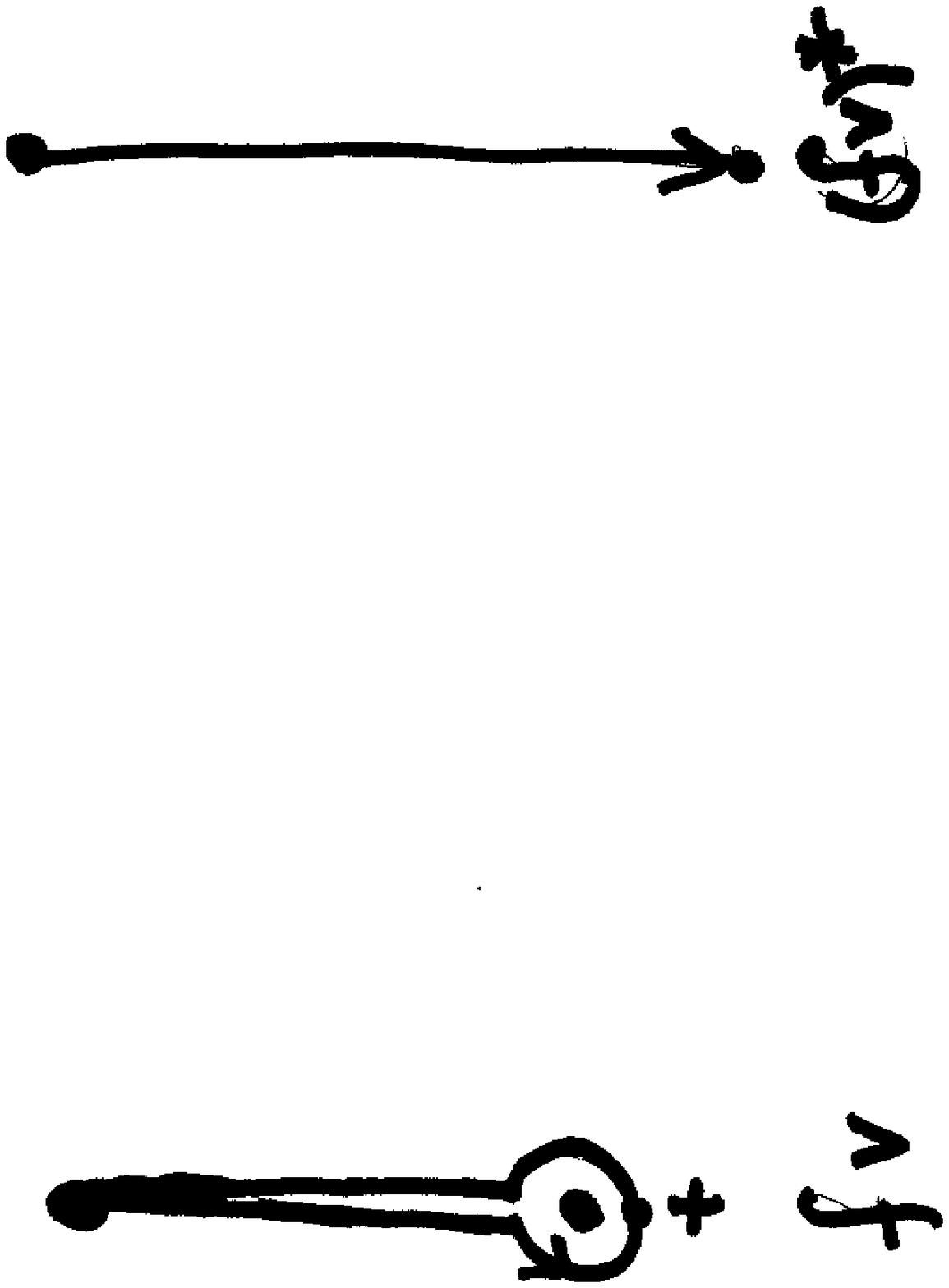}
\vskip 6cm \botcaption{Figure 2}
   Chains and quantum group cf. [8,11,29]. Phase is chosen to be zero
at the point marked with $+$. \text{\footnotemark} 
\endcaption
\endinsert
\vskip 6cm

\midinsert\vskip 0.5cm
\includegraphics{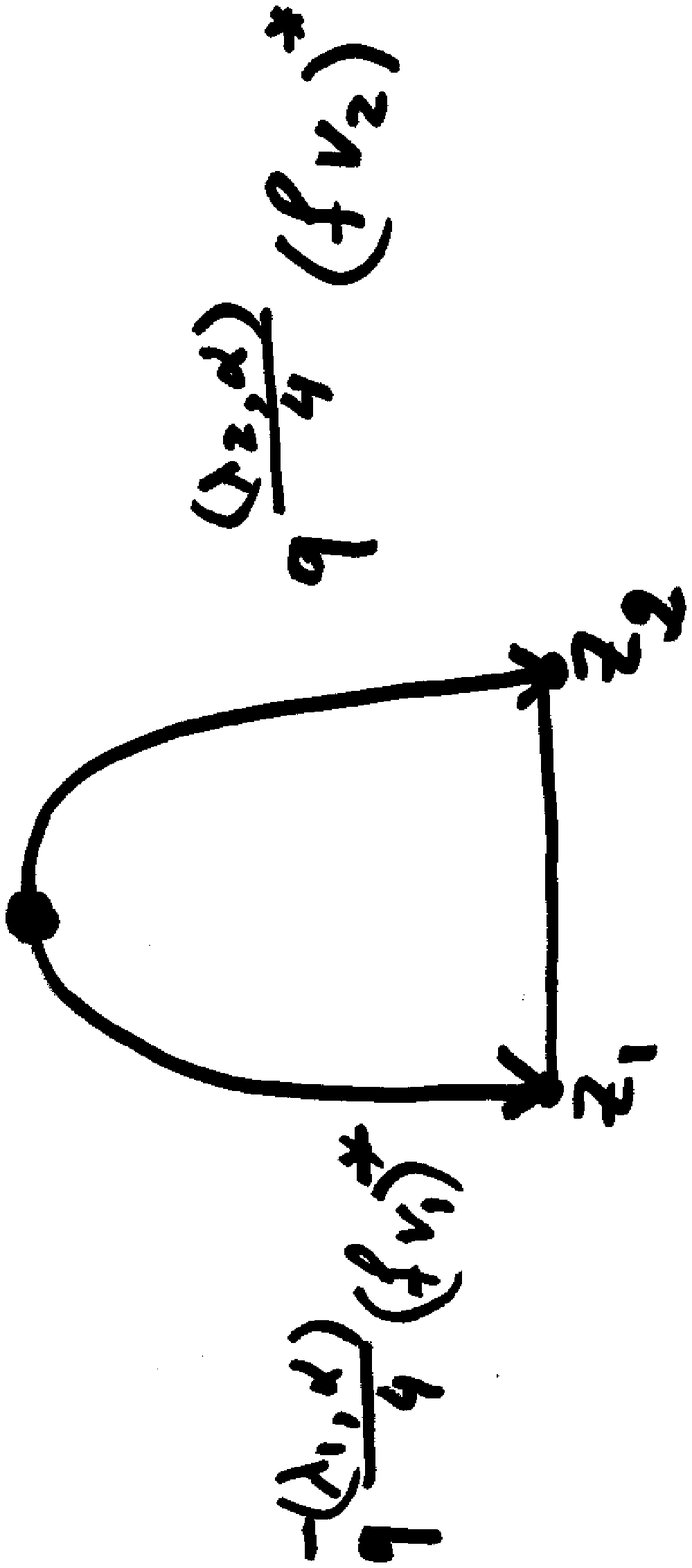}
\vskip 6cm \botcaption{Figure 3}
   Decomposition with the help of quantum group cf. [8,11,29].
\endcaption
\endinsert

\vskip 6cm
\footnotetext{The choice of comultiplication is dictated by the choice
that the phase is equal to zero at the point marked with $+$.}

\subhead {5.6. Case of root system of type $A_n$ }\endsubhead

\proclaim{Theorem 5.7}
Let $v_0, v_1,\ldots ,v_{n+1}$ be the highest weight vectors of
irreducible highest weight
 modules with highest  weights $\Lambda(0),
\Lambda(1)=\Lambda(2)=\ldots=\Lambda({n+1})=\Lambda$ correspondingly.
Here $\Lambda(0), \Lambda$ are as in section 5.1. Then
cycle  $\pmb\Delta $ is encoded as:

1.
$$
\multline
v_{\Delta}=
\sum_{w \in S_{n+1}} 
(-1)^{ l(w)} q^{{1\over 4} [{n(n+1) \over 2} - 2 l(w)]}
v_0^{*} \otimes (f_{w(1)-1} \ldots f_2 f_1 v_1)^{*} \otimes \ldots \\
\otimes (f_{w(i)-1} \ldots f_2 f_1 v_2)^{*} \otimes 
\ldots \\ \otimes(f_{w(n+1)-1} \ldots f_2 f_1 v_{n+1})^{*}
\endmultline
$$

2.
$$
\multline
v_{\Delta}=  \sum_{w \in S_{n+1}} \frac{(-1)^{ l(w)} q^{{1\over 4}[{n(n+1)\over 2} -2
l(w)]}}
{(q^{1\over 2} - q^{-{1\over 2}})^{n(n+1)\over 2}} \\ \times
v_0 \otimes f_{w(1)-1}\ldots  f_{2} f_1  v_1 \otimes
\ldots \\ 
\otimes f_{w(i)-1}\ldots  f_{2} f_1 v_i 
\otimes \ldots \\
\otimes  f_{w(n+1)-1}\ldots  f_{2} f_1 v_{n+1}
\endmultline
$$
In particular, $e_i v_{\Delta} =0$ for $i=0,1,\ldots n$,
this implies that cycle $\pmb\Delta$ defines  correctly an element of
$H_{*!,m}(U,S^{*})^{-}$ and the use of  the word `cycle' is justified. 
\endproclaim
\medskip
To prove the theorem we repeatedly use figures 3 and 2.
Also for the above theorem  the theorem 2.6 and remark 2.7 of [17] are
helpful : the number of arrows which are `to the left' is equal to the
length of the corresponding element of the Weyl group (also reproduced
below in theorem 5.3). In fact using
elementary decomposition fig. 3 one gets that each arrow to the right
brings the factor $q^{1\over 4}$ ,while each arrow to the left brings the
factor $(-1) q^{-1\over 4}$. The number of variables $t_i^{(j)}$ is equal
to $\frac{n(n+1)}{2}$. The key point is that
all the `wrong' diagrams (cf. fig. 5) cancel each other
because of the phase argument cf. fig. 6.
 Where  by the `wrong' diagrams we mean
the diagrams in which there are two arrows with the same target.

Vice versa, this theorem might be considered as quantum group
explanation of theorem 2.6 of [17] (numbers of arrows which are to the
left $=l(w)$).

\smallskip
\subhead{5.8. Diagrams} \endsubhead
 The notion of a diagram was introduced in [38]
in the context of trigonometric Knizhnik-Zamolodchikov equation, and 
later similar
notion was introduced in [53] in the context of
 multidimensional determinants and discriminants. We will use the
notion of a
 diagram in the
form of [53], in particular, we borrow the very convenient graphical
notation (see figs. 1,2,3 of [53]).

\medskip 
Fix some positive integer $n$. Consider the set of $\frac{n(n+1)}{2}$
points, indexed by pairs of integers
 $\{(i,j)|$ $i=1, \ldots,j$, $j=1,\ldots, n \}$. 
It is helpful to organize the points in the form
of a pattern, so that points are divided in $n$ rows, $j$th row is formed
by points $\{(i,j)| \; i=1, \ldots,j\} $; point $(i,j)$ is located under and
between points $(i,j+1)$ and $(i+1, j+1)$ (fig. 4a).

\midinsert\vskip 4cm
\includegraphics{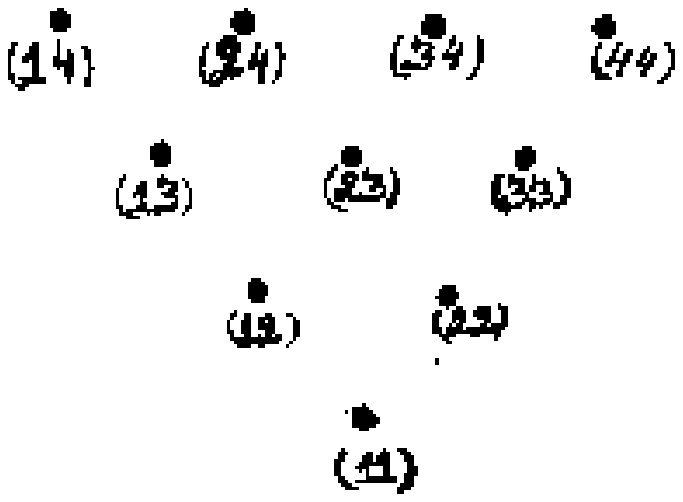}
\botcaption{Figure 4a}
$n=4$.\quad Points $ (i,j)$ organized in a pattern;
$j$ is the number of the row, $i$ is the number in the row
\endcaption
\endinsert

Now mark with a cross one point in each row. Let $\{(i_j,j)\}$ be the
subset of marked points (fig. 4b).
\vskip 6cm
\midinsert\vskip 7cm
\includegraphics{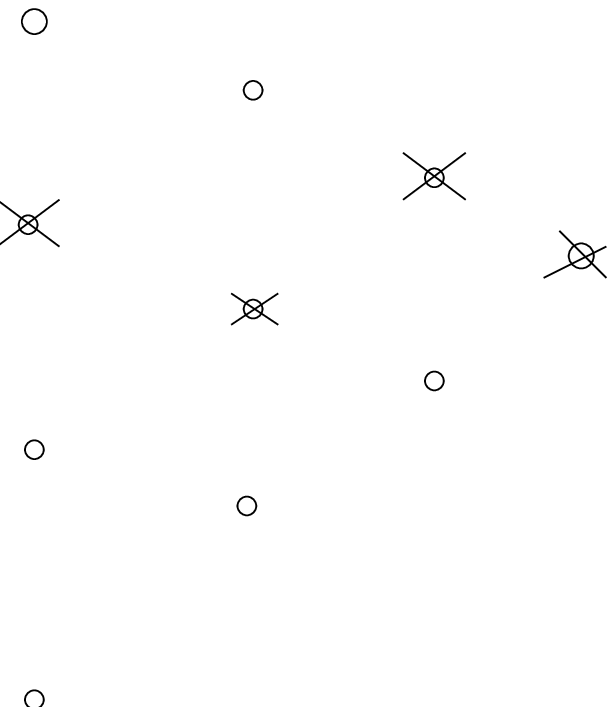}
\botcaption{Figure 4b}
One point is marked in each row.
\endcaption
\endinsert

Finally, draw an arrow for each point $(i,j)$ with the source in this
point $(i,j)$ and target $tar(i,j)$ in the next $j+1$th row defined as:
$$
tar(i,j)=
\cases
(i,j+1), &\text{if }i<i_{j+1}\\
(i+1,j+1), &\text{if }i \geq i_{j+1} \text{\footnotemark}
\endcases
$$
\footnotetext{Recall that $\{(i_j,\;j) \} $ is the set of marked points and
 $ (i_{j+1},\; j+1)$
is the only marked point in $j+1$th row.}
If $tar(i,j)=(i,j+1)$, then the arrow is called to the {\bf{left}},
if $tar(i,j)=(i+1,j+1)$, then the arrow is called to the {\bf{right}} .

Note: neither arrow has a marked point as its target.

In this way one obtains fig 4c.

\midinsert\vskip 6cm
\includegraphics{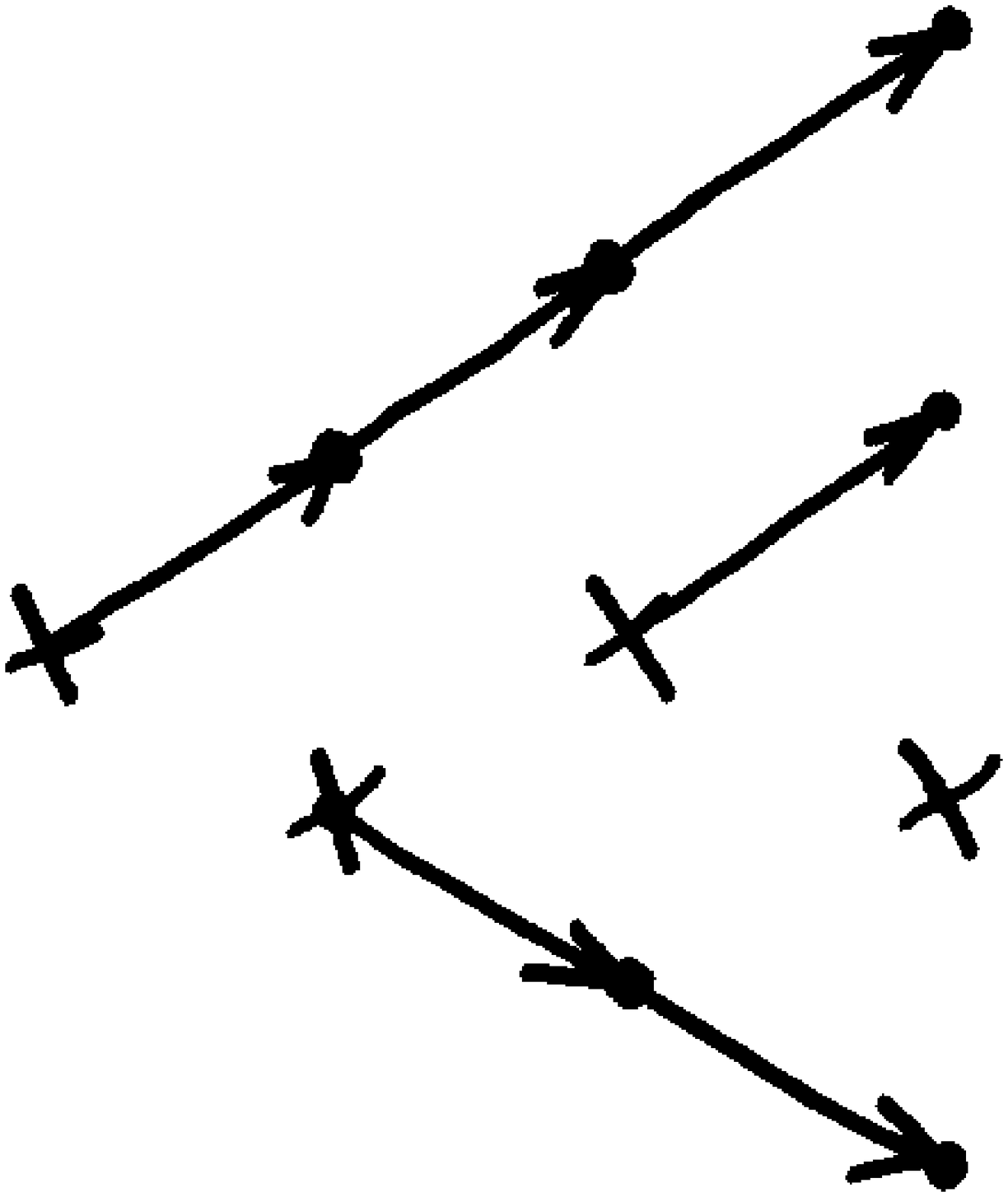}
\botcaption{Figure 4c}
Example of a diagram.
\endcaption
\endinsert

\definition{Definition 5.9}
A triple $$ (\{(i,\;j)\},\quad \{(i_j,\;j)\},\quad tar) $$
 consisting of:

set of points
 $\{(i,\;j)| \quad i=1,\ldots,j\; ; \; j=1, \ldots,n\} $,

set of marked points $\{ (i_j,\;j) \; | \quad j=1,\ldots,n \} $ 

and function 
$\quad tar \quad $ defined above will be called {\bf a diagram.} \enddefinition

\remark {Remark 5.10} One can see that a diagram is determined by the set of
marked points.
\endremark

\definition{Definition 5.11} We describe the correspondence between
diagrams and elements of Symmetric group as follows.
Consider a diagram as an oriented graph and forget orientation.
For $i=1, \ldots,n$ define $ w(i)$ as the number of points in the connected
component of the point $(i,n)$.
\enddefinition
\smallskip

Symmetric group $S_n$ has standard generators
 $\sigma_1,\;\sigma_2,\ldots, \sigma_{n-1}$, where
 $\sigma_i$ permutes $i $ and $ i+1$.

\definition{Definition 5.12} The length $l(w)$ of an element 
$w \in{S_n}$ is 
the minimal
integer $p \ge 0$, s.t. $w$ admits a presentation 

$$ w = \sigma_{i_1}\sigma_{i_2} \ldots \sigma_{i_p}  \;.$$
Any presentation of $w$ as a product of $p=l(w)$ generators is called
a reduced presentation.

\enddefinition

\proclaim {Theorem 5.13}
 Let a diagram
$$
(\{(i,j)| \quad i=1, \ldots,j,\; j=1, \ldots,n\}, \quad \{(i_j,j) | \quad
j=1, \ldots,n \}
 \quad, \; tar)
$$
corresponds to an element $w\in{ S_n}$. Then the length $ l(w)$ of an
element $w$ is equal to:
$$
l(w)=\sum_{j=1}^n{(i_j -1)}
$$
In other words , $l(w)$ is equal to the number of arrows which are to
the left in the diagram corresponding to the element $w$.
cf. theorem 2.6 of [17].
\endproclaim

\midinsert\vskip 0.5cm
\includegraphics{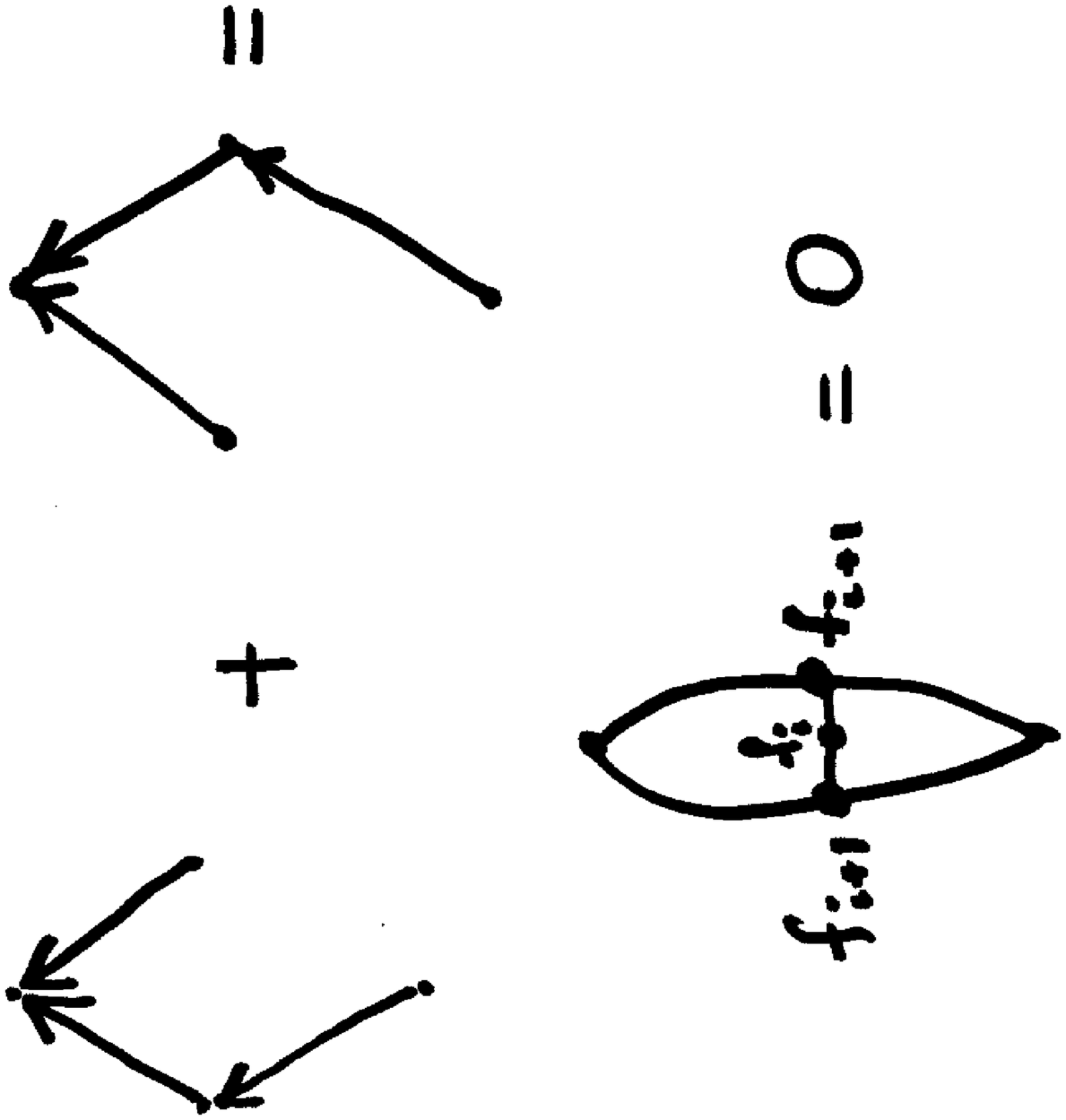}
\vskip 8cm \botcaption{Figure 5}
 Wrong diagrams cancel each other (quite in the same way this works
in Harish-Chandra decomposition [18])
\endcaption
\endinsert
\vskip 6cm

\midinsert\vskip 15cm
\includegraphics{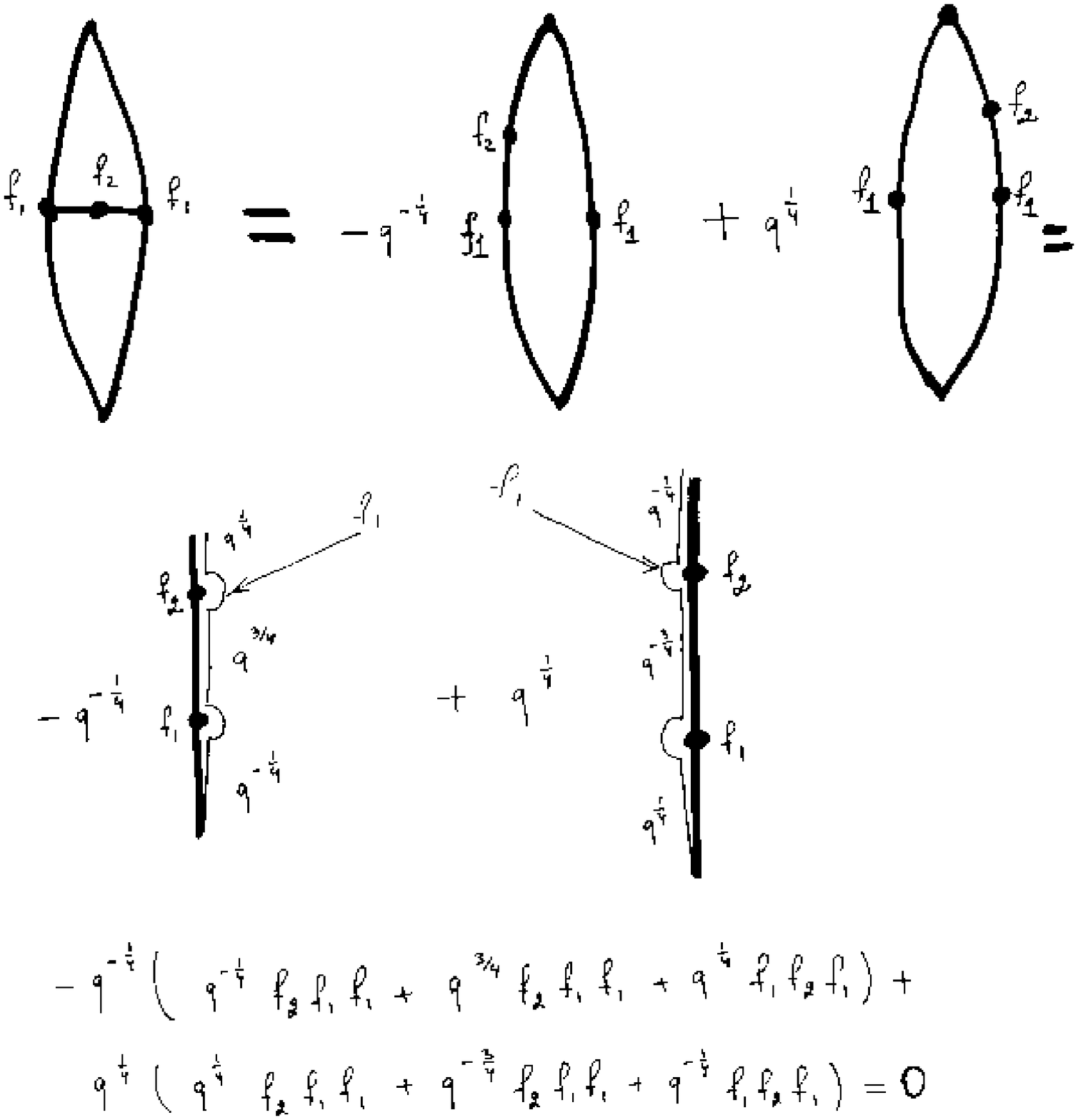}
 \botcaption{Figure 6}
     Phase argument: decomposition into ordered chains-
     total sum is identical zero chain
\endcaption
\endinsert
\vskip 6cm

\smallskip
\proclaim{Theorem 5.14}  Let $\Lambda(1)=\Lambda(2)=\Lambda$ as in
theorem 5.7. Also, let $v_1,v_2$ be the highest weight vectors of
Verma modules $M(\Lambda(1)), M(\Lambda(2))$, correspondingly.
Let $M(\Lambda(1))^{*}, M(\Lambda(2))^{*}$ be dual Verma modules.  
Let $PR$ be the braiding: 
$$
PR \:
M(\Lambda(1))^{*} \otimes M(\Lambda(2))^{*} \mapsto M(\Lambda(2))^{*}\otimes
M(\Lambda(1))^{*}
.
$$
Then for 
$i>j$ one has:
$$
\multline
PR \; ((f_i f_{i-1}\ldots f_1 v_1)^{*} \otimes
(f_j f_{j-1} \ldots f_1 v_2)^{*}) = \\
q^{1\over 2} (f_j f_{j-1}\ldots f_1 v_2)^{*} \otimes
(f_i f_{i-1} \ldots f_1 v_1)^{*} + \\
q^{1\over 2}( q^{1\over 2} - q^{-{1\over 2}})
(f_i f_{i-1} \ldots f_1 v_2)^{*} \otimes
(f_j f_{j-1} \ldots f_1 v_1)^{*}
\endmultline
$$

For $i<j$ one has:
$$
\multline
PR \; ( (f_i f_{i-1} \ldots f_1 v_1)^{*}\otimes
(f_j f_{j-1} \ldots f_1 v_2)^{*})= \\
q^{1\over 2} (f_j f_{j-1} \ldots f_1 v_2)^{*}\otimes
(f_i f_{i-1} \ldots f_1 v_1)^{*}
\endmultline
$$

\endproclaim
\smallskip
Proposition is easily proved by contour manipulations, see fig. 7
as elementary, but typical example. All calculations are up to the
kernel of contravariant form $S$.

\midinsert\vskip 0.5cm
\includegraphics{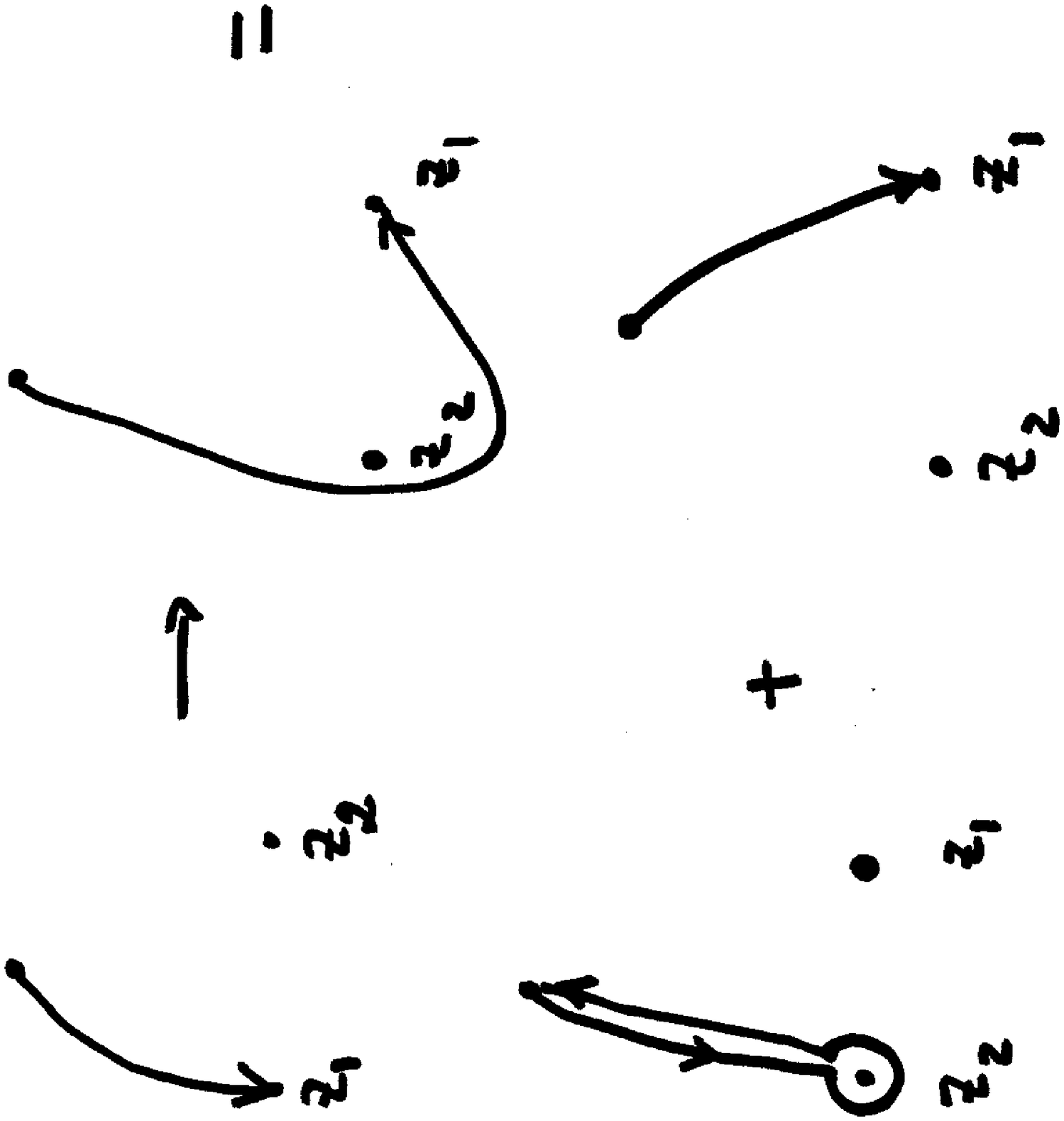}
\vskip 6cm \botcaption{Figure 7}
      Result of the braiding
\endcaption
\endinsert

\vskip 3cm
 
\proclaim{Corollary 5.15} Suppose $ z_1(t), z_2(t), \ldots, z_{n+1}(t)$, $t \in
[0,1]$ are closed loops on a complex plane, i.e. $z_1(0)=z_1(1), 
z_2(0)=z_2(1), \ldots , z_{n+1}(0)=z_{n+1}(1)$, such that $z_i(t) \ne z_j(t)$
for $i \ne j$. Let also $Re(z_i(t)) >0$ for each $i=1,\ldots ,n +1$.
Then the homological class of the cycle $\pmb\Delta$ is
preserved under 
the
monodromy along paths $z_i(t)$.
\endproclaim
\smallskip
In fact such monodromy can be produced as composition of even number
of elementary braidings as in theorem 5.14, each of them gives the
factor $-1$. 
Braiding of $z_1$ with
$z_0=0$ is forbidden by hypotheses.

\subhead {5.16. Using R-matrix for tensor product of  
vector representations} \endsubhead
R-matrix for the tensor product of vector representations of
$sl(n+1)$
 reads
as:
$$
R=q^{-1/2(n+1)} \Big\{ \sum_{i \ne j} E_{ii} \otimes E_{jj}
+q^{1\over 2} \sum_{i}E_{ii}\otimes E_{ii} +
(q^{1\over 2}  - q^{-{1\over 2}}) \sum_{i < j} E_{ij} \otimes E_{ji} \Big\} 
$$
cf. [26, 27]. Vector representation of $sl(n+1)$ has natural basis 
$e_1, e_2, \ldots , e_{n+1}$
with highest weight vector 
$e_1$,  
$E_{ij}$
are matrix units:

$$
E_{ij}e_k= \delta_{kj}e_i ,
$$
where $\delta_{kj}$ is a Kronecker's delta.
In our case the interesting part of R-matrix can be easily obtained from the above one
by multiplying by $q^{{n+2}\over {2(n+1)}}$:

$$
\hat R = {q^{1\over 2}}{ \Big\{ \sum_{i \ne j}{ E_{ii} \otimes E_{jj}}+
{q^{1 \over 2}} \sum_{i}{ E_{ii}\otimes E_{ii}} +
(q^{1\over 2} - q^{-{1\over 2}}) \sum_{i < j}{ E_{ij} \otimes E_{ji} } \Big\}}.
$$

Then one immediately verifies that:

$$
\hat R ( {e_k \otimes e_l})= q^{1\over 2}{ e_k \otimes e_l} \quad \text{if}
\quad 
k<l
$$

and 
$$
 \hat R (e_k \otimes e_l) = q^{1\over 2} e_k \otimes e_l +
q^{1\over 2} ( q^{1\over 2} - q^{-{1\over 2}}) e_l \otimes e_k
\quad if \quad  l < k
$$

Also,
$$ \hat R (e_k \otimes e_k) = q e_k \otimes e_k$$

Let $P$ be the  transposition, i.e. 
$$P e_i \otimes e_j = e_j \otimes e_i$$

Then one immediately obtains that for $k < l$
$$
{P \hat R} (q^{1\over 4} e_k \otimes e_l - q^{- {1\over 4}} e_l \otimes
e_k)=
(-1)( q^{1\over 4} e_k \otimes e_l - q^{-{1\over 4}} e_l \otimes e_k) 
$$
i.e. 
$$
q^{1\over 4} e_k \otimes e_l- q^{-{1\over 4}} e_l \otimes e_k
$$
is an eigenvector of $P \hat R$ with eigenvalue $-1$ !

\smallskip
\subhead{5.17 Tensor product with vector representation} \endsubhead
As it is shown by Lusztig and Rosso [46,47] representation theory of
$U_q( g)$ for generic values of $q$ is the same as in the classical
case $q=1$.
The tensor product of finite-dimensional representations of $sl(n+1)$ is
governed by the Littlewood-Richardson rule, cf. [45] , which is in turn
a consequence of the theory of characters and symmetric functions.
The rule is essentially simple for the tensor product with
vector representation, namely, we should add one box to the Young
diagram so that in result we obtain again a Young diagram and if the column
with $n+1$ boxes appears it should be removed (as corresponding to trivial 
representation). The  tensor product with vector 
representation  is multiplicity free.
For example, the tensor product of two vector representations
decomposes as:
$$
L(\eta_1)\otimes L(\eta_1) = L(2 \eta_1) \oplus
L(\eta_2)
$$
where $L(2 \eta_1)$, $L(\eta_2)$ denotes 
the analog of the symmetric tensor, antisymmetric tensor, see fig. 8.

\midinsert\vskip 0.5cm
\includegraphics{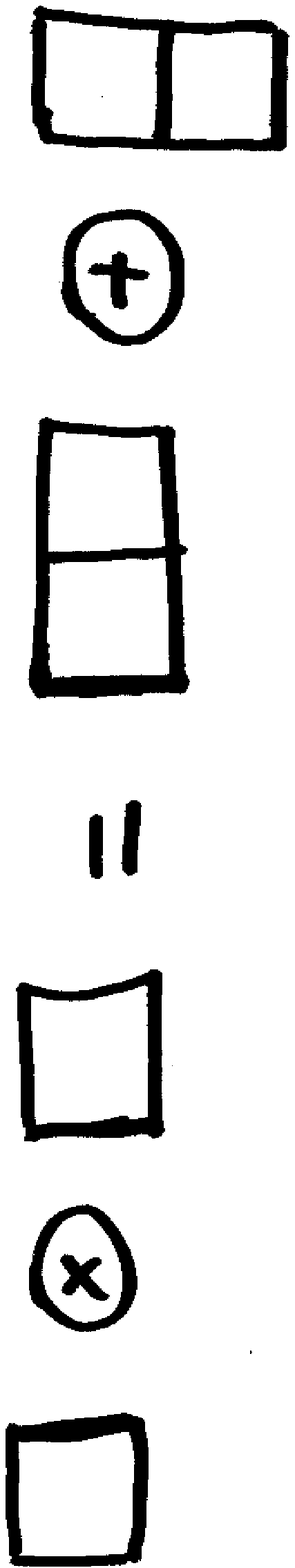}
\vskip 3cm \botcaption{Figure 8}
   Tensor product of two vector  representations, 
          Littlewood-Richardson rule
\endcaption
\endinsert

Let $P_{2 \eta_1}, P_{\eta_2}$ denotes  the
corresponding  projectors.
The PR-matrix is given as follows:
$$
PR = q^{-\frac {1}{2(n+1)}}(q^{1\over 2} P_{2 \eta_1}- 
q^{-1\over 2}P_{\eta_2}) .$$

So one can see that on antisymmetric tensors it acts as 
$-q^{-\frac{n+2}{2(n+1)}}$.
Multiplying by $q^{\frac{n+2}{2(n+1)}}$ as in previous section we get
the coefficient of $P_{\eta_2}$ is
 $-1$ .
\smallskip
Now consider the product of finite-dimensional representation of
weight $\lambda$ with vector  representation:
$L(\lambda) \otimes L(\eta_1)$. And assume that we add the box
to the $s$-th row of a Young diagram corresponding to $\lambda$. In
particular , we assume that this is a correct operation.
This means that we have variables of integration 
$\{t_i | i=1, \ldots, s-1\}$ and the integral:

$$
\int \prod_{i=1}^{s-1} t_i^{\frac{(\lambda, -\alpha_i)}{\varkappa}}
\prod_{j < i} (t_i -t_j)^{\frac{(-\alpha_i, -\alpha_j)}{\varkappa}}  
\prod_{i=1}^{s-1} (z - t_i)^{\frac{(\eta_1,-\alpha_i)}{\varkappa}}
\frac {d t_1}{t_1} \frac {dt_2} {t_2} \ldots \frac {dt_{s-1}}{t_{s-1}}
$$
The natural domain of integration for asymptotic solution is:
$$
0 \le t_{s-1} \le t_{s-2} \le \ldots \le t_1 \le z.
$$

\midinsert\vskip 7cm
\includegraphics{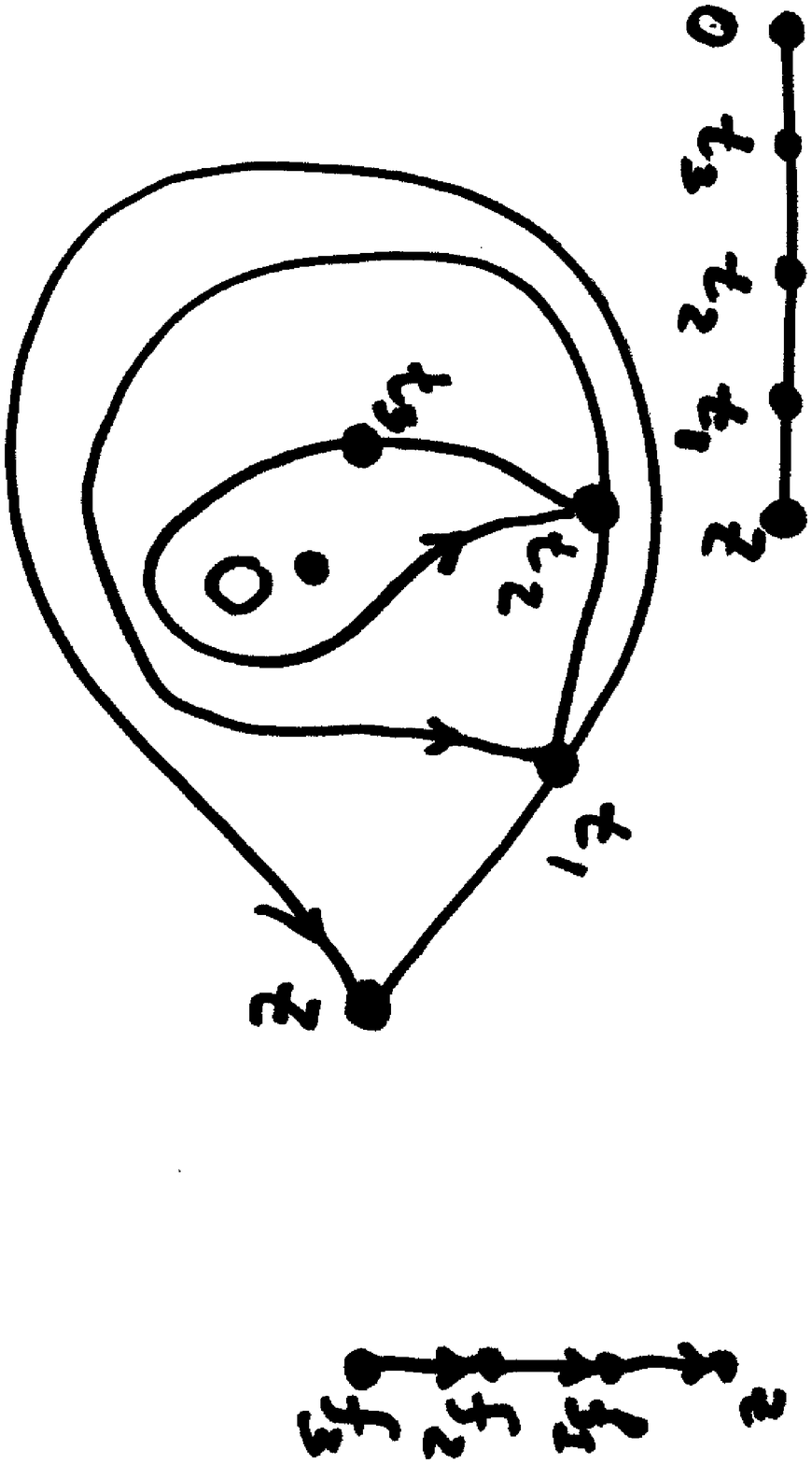}
\vskip 2cm \botcaption{Figure 9}
  Cycle for integration for screened  vertex fundamental
operator: contour  for $t_i$ starts and ends at $t_{i-1}$,$i=2,3,\ldots$;
contour for $t_1$ starts and ends at $z_1$. In particular, all
'internal' contours are  movable. 
\endcaption
\endinsert

The leading asymptotic is equal to
$$z^{\frac{(\lambda,-\sum_{i=1}^{s-1} \alpha_i) - (s-1)}{\varkappa}}.$$
The integral is easily taken using the formula 
for Euler's beta function:
$$
\int_0^1 t^{\alpha-1} (1-t)^{\beta -1} d t =
\frac {\Gamma(\alpha) \Gamma(\beta)}{\Gamma(\alpha + \beta)}
$$ 
and the leading asymptotic coefficient  is equal to:
$$
\prod_{p=1}^{s-1}
\frac{\Gamma \Big( \frac{(\lambda,-\sum_{i=s-p}^{s-1}
       \alpha_i)} {\varkappa}  -\frac{(p-1)}{ \varkappa}\Big) 
       \Gamma(1-\frac{1}{\varkappa})}
{\Gamma  \Big(\frac {(\lambda,-\sum_{i=s-p}^{s-1} \alpha_i)}{\varkappa} +1
-\frac{p}
{\varkappa}\Big)}
$$
The above  cycle for integration can be encoded as  
singular vector of
tensor
product of two irreducible modules  with highest weights 
$\lambda$ and $\eta_1$
as follows:
$$
\align
&(-1)^{s-1}q^{{\frac{-1}{4}(s-2)+ \frac{1}{4}(\lambda, \sum_{i=1}^{s-1} -\alpha_i )}} (f_{s-1} \ldots f_1 v_0)^{*} \otimes v_1^{*} +\\
&q^{\frac{1}{4}(s-1) } v_0^{*} \otimes (f_{s-1} \ldots f_1 v_1)^{*}+\\ 
&\sum_{p=1}^{s-2}  
(-1)^{s-p-1} q^{{\frac{-1}{4}(s-p-2)}+{\frac
{1}{4}p}+{\frac{1}{4}(\lambda,\sum_{p+1}^{s-1} - \alpha_i)}}
(f_{s-1} f_{s-2} \ldots f_{p+1}  v_0)^{*} \otimes 
(f_p f_{p-1} \ldots f_1 v_1)^{*}
\endalign 
$$

With the help of comultiplication $\Delta$ this might be iterated 
(in principle) and
cycles for asymptotic solutions from ref. [17] can be encoded as
 singular vectors  of the tensor product of irreducible
highest weight modules over quantum group:
 
$$L(\Lambda(0)) \otimes L(\Lambda(1)) \otimes \ldots \otimes
L(\Lambda(n+1)).$$
In the case of generic $\Lambda(0)$ one can use
Bernstein-Gelfand-Gelfand results on the
category $\Cal O$ cf. [78], see also Kostant [77] ( and the theory of
Knizhnik-Zamolodchikov equations, parameter $q$ is assumed to be generic).
In particular for generic $\Lambda(0)$ tensor product $L(\Lambda(0))
\otimes L(\eta_1)$  decomposes as:
$$
L(\Lambda(0))
\otimes L(\eta_1)=L(\Lambda(0) + h_1) \oplus L(\Lambda(0)+h_2) 
\oplus \ldots \oplus L(\Lambda(0)+ h_{n+1})
$$
The space of singular vectors of weight $\Lambda(0)$
of the tensor product
$$L(\Lambda(0)) \otimes L(\eta_1) \otimes \ldots \otimes L(\eta_1)$$
with $n+1$ factors $L(\eta_1)$ is $(n+1)!$-dimensional.

%Also the cycle 
% can be decomposed into elementary chains of Felder's
% %type [10]. Though elementary chains do not  belong  to the
%homology group, their linear combination does.
%
%% 
\smallskip

\remark{Remark 5.18} The fact that the cycle $\pmb\Delta$ corresponds to
q-antisymmetric tensors is actually clear. Consider  tensor
product of $n+1$ vector representations.
So  there is the only way to get the
variables of integration as needed: namely, to {\it add box under box,
so that the Young diagram will be the column of $n+1$ boxes}, which corresponds
to antisymmetric tensors ,see fig. 10. 
One could also find different contours for integration 
giving the same homological class as cycle $\pmb\Delta$ (fig. 11), but cycle
$\pmb\Delta$ is especially
convenient for obtaining Harish-Chandra decomposition, it
has geometric origin in harmonic analysis and there is certain parallellism
with Gelfand-Tsetlin patterns (see also [73, 40],[74]).
\endremark
\smallskip

\vskip 6cm
\midinsert\vskip 2cm
\includegraphics{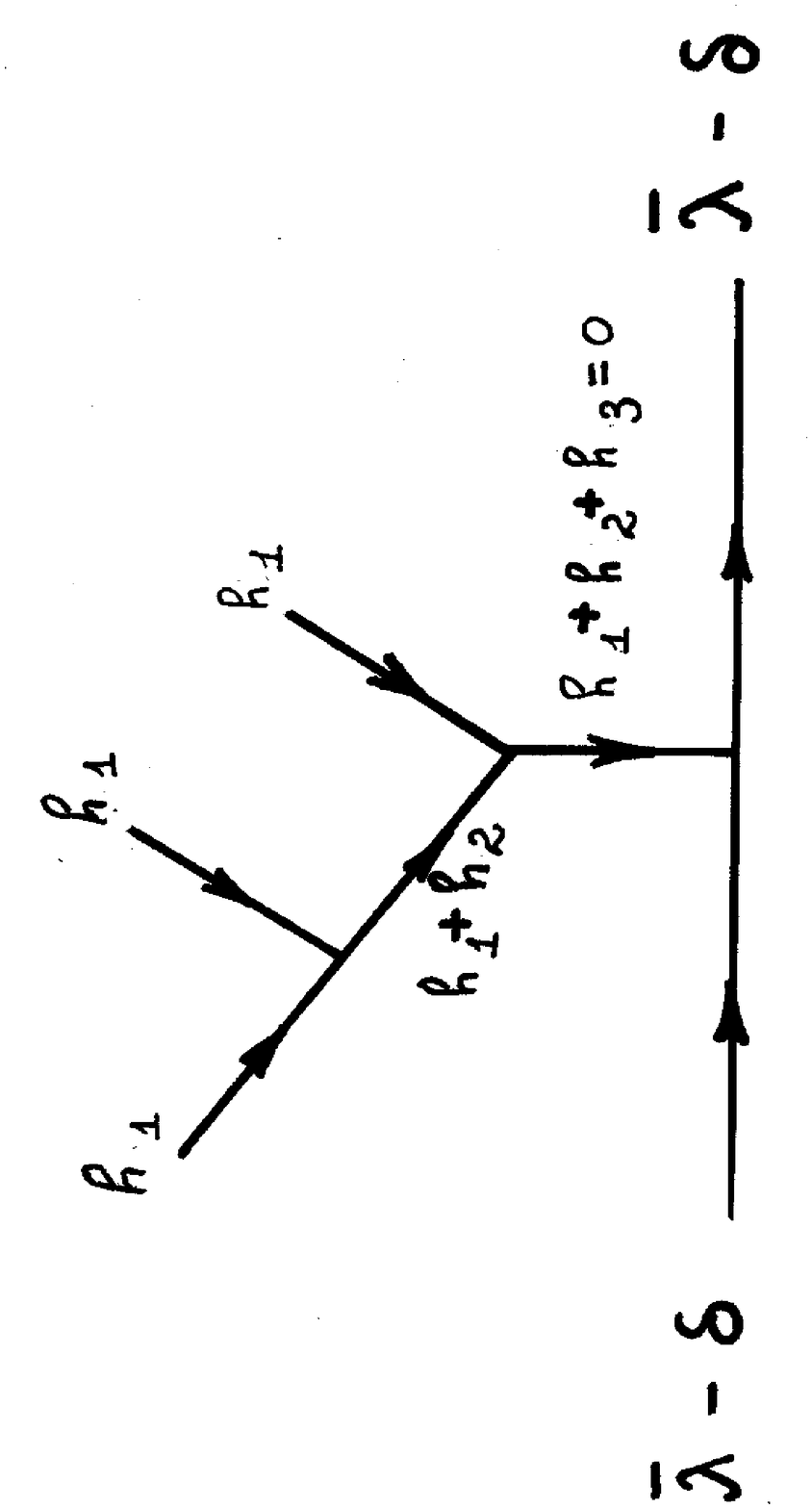}
\vskip 6cm \botcaption{Figure 10}
 Conformal block  for zonal spherical
function ( $A_2$ case). 
\endcaption
\endinsert

\vskip 15cm 
\midinsert\vskip 7cm
\includegraphics{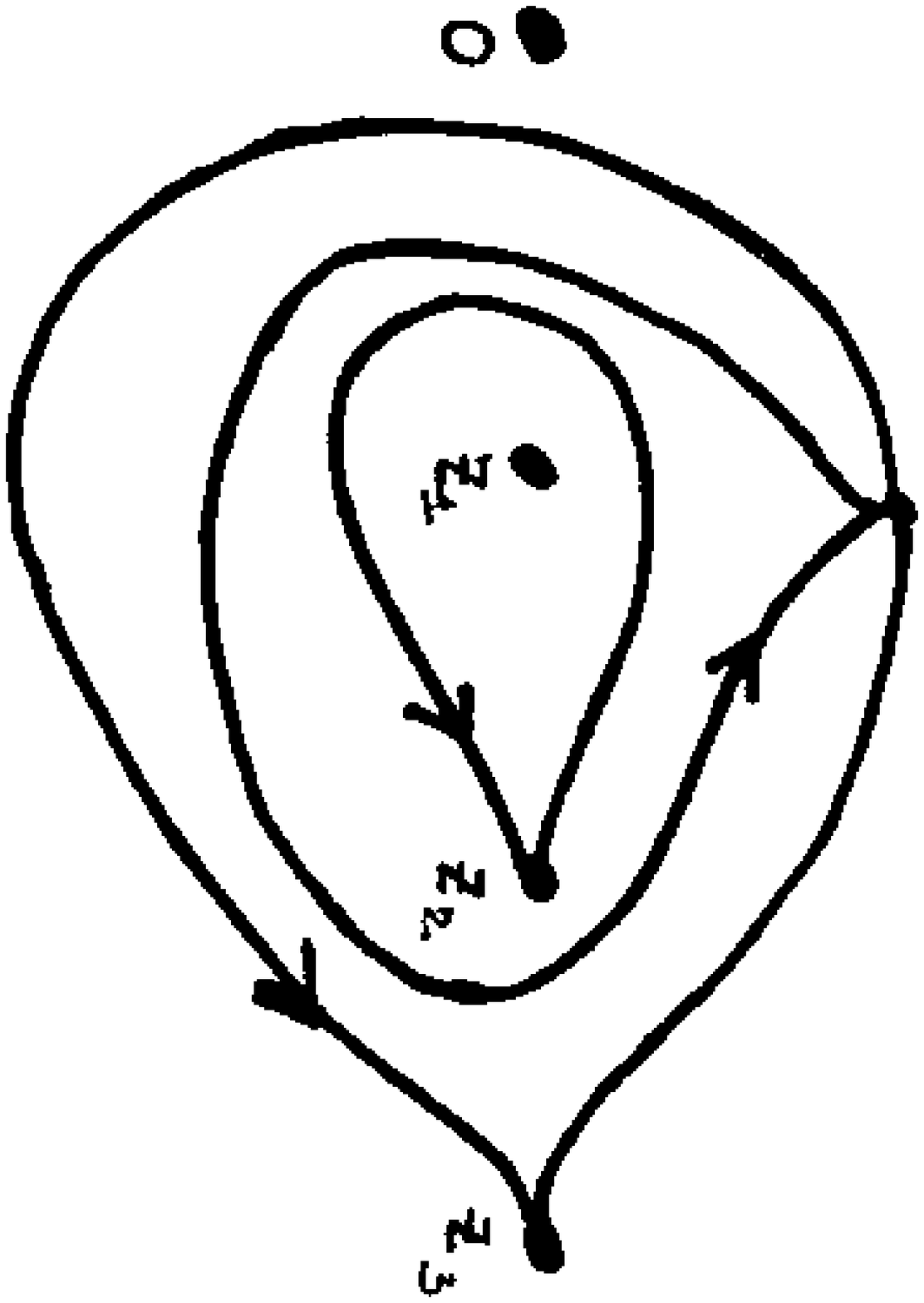}
 \botcaption{Figure 11 }
  Another contour for integration for zonal spherical function. 
\endcaption
\endinsert

\subhead{Concluding remarks }\endsubhead
 The distinguished  cycle
$\pmb\Delta$ serves as a contour for integration for zonal spherical
function of type $A_n$. It goes back to  the classical calculation
of Gelfand and Naimark of zonal spherical function for $SL(n,\Bbb C)$,
originates in the so-called elliptic coordinates and provides a
materialization of the flag manifold. Zonal spherical function is a
particular conformal block. 

\subhead Acknowledgments \endsubhead
We are grateful to  I.~Gelfand for the suggestion to use the
technique of ref. [23] in application to Heckman-Opdam hypergeometric
functions, to S.~Lukyanov for discussions devoted to conformal field
theory and $W$-algebras, to V.~Schechtman for discussions
concerning refs. [8,12,57], to V.~Brazhnikov, M.~Braverman for
stimulating discussions.

\enddocument